\documentclass{article}
\usepackage[utf8]{inputenc}
\usepackage{amsmath}
\usepackage{amsfonts}
\usepackage{algorithm}
\usepackage{algpseudocode}
\usepackage{graphicx}
\usepackage{hyperref}
\usepackage{makecell}
\usepackage{caption}
\usepackage{subcaption}

\newtheorem{theorem}{Theorem}

\newtheorem{conjecture}{Conjecture}

\begin{document}

\title{Optimizing for the Rupert property}
\markright{Optimizing for the Rupert property}
\author{Albin Fredriksson\footnote{albin.fredriksson@gmail.com}}

\maketitle

\begin{abstract}
\noindent
A polyhedron is Rupert if it is possible to cut a hole in it and thread an identical polyhedron through the hole. It is known that all 5 Platonic solids, 10 of the 13 Archimedean solids, 9 of the 13 Catalan solids, and 82 of the 92 Johnson solids are Rupert. Here, a nonlinear optimization method is devised that is able to validate the previously known results in seconds. It is also used to show that 2 additional Catalan solids---the triakis tetrahedron and the pentagonal icositetrahedron---and 5 additional Johnson solids are Rupert.
\end{abstract}

\noindent
\section*{Introduction}
A convex polyhedron\footnote{Polyhedra discussed in this paper are always assumed to be convex.} is said to have the \emph{Rupert property} (or to ``be Rupert'') if it is possible to cut a hole in it and thread another identical polyhedron through the hole. At first glance, it can appear quite counterintuitive that this should be possible, but many polyhedra have been shown to possess the property, and no polyhedron has been shown to constitute a counterexample.

\begin{figure}[htb]
  \centering
    \begin{subfigure}[t]{0.32\textwidth}
     \centering
     \includegraphics[width=\linewidth]{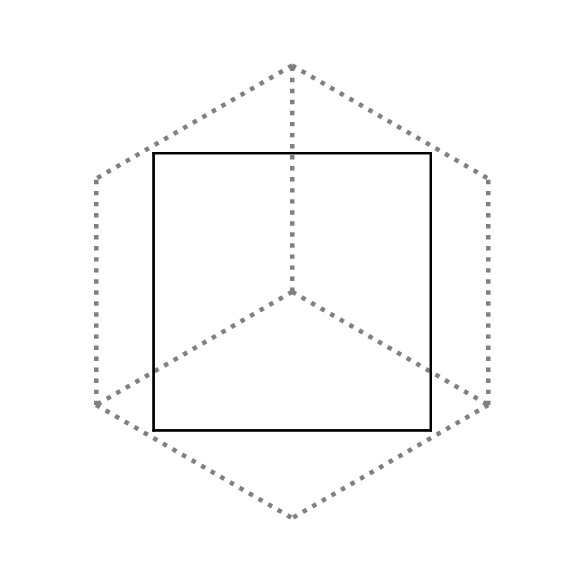}
     \caption{The passage through the cube seen from a corner.}
     \label{fig-cube1}
    \end{subfigure}
    \hspace{0.5cm}
    \begin{subfigure}[t]{0.32\textwidth}
     \centering
     \includegraphics[width=\linewidth]{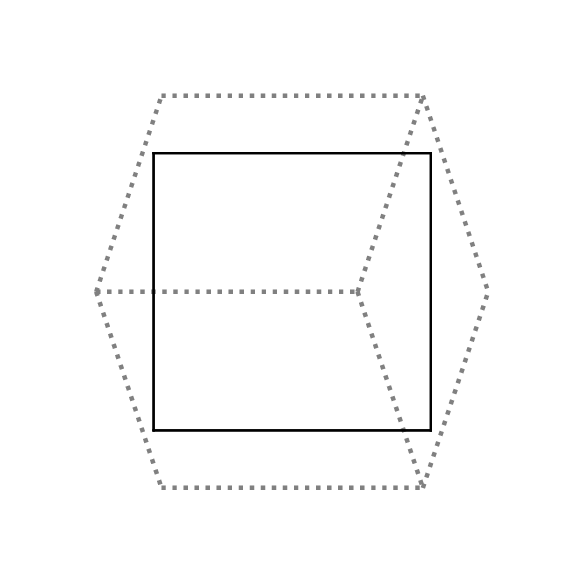}
     \caption{The passage through the cube oriented to provide the largest possible margin to the passing cube.}
     \label{fig-cube-opt}
    \end{subfigure}
  \caption{Rupert passages through the unit cube. Black solid lines represent cubes passing through those represented by gray dotted lines.}
  \label{fig-cube}
\end{figure}

The history of the Rupert property dates back to the 17th century, when Prince Rupert of Rhine wagered that a hole can be cut in the unit cube such that another unit cube can be threaded through it and won~\cite{schreck}. To see that such a threading is possible, one can look at the unit cube from one of its corners as in Figure~\ref{fig-cube1}---the unit square actually fits with a margin. A slight tilting of the cube yields an even larger margin, see Figure~\ref{fig-cube-opt}. The latter passage was found by Pieter Nieuwland in the 18th century. The largest scale factor for which a scaled copy of a polyhedron can be threaded through an unscaled copy is therefore called the \emph{Nieuwland number} of the polyhedron.

In 1968, Scriba~\cite{scriba} showed that the tetrahedron and the octahedron are Rupert. Many additional results have been obtained in recent years. In 2017, Jerrard et al.~\cite{jerrard} showed that the remaining Platonic solids are also Rupert. In 2018, Chai et al.~\cite{chai} showed that 8 of the 13 Archimedean solids are Rupert, and in 2019, the same was shown for a ninth Archimedean solid~\cite{hoffmann, lavau}. A higher-dimensional result was given by Huber et al.~\cite{huber}, who showed that the $n$-cube is Rupert.

Even more recently, Steininger and Yurkevich~\cite{steininger} used an algorithm based on checking whether random orientations of polyhedra permit Rupert passages to show the property for a tenth Archimedean solid, the icosidodecahedron, as well as 9 of the 13 Catalan solids (the duals of the Archimedean solids), and 82 of the 92 Johnson solids. They improved on the lower bounds of the Nieuwland numbers of their solutions by an optimization procedure in which they iteratively applied random perturbations to a given solution and accepted a perturbed solution if it led to an improvement. Moreover, they devised a deterministic algorithm that can theoretically determine whether a given polyhedron is Rupert or not. It requires finding or showing the nonexistence of solutions to many simultaneous polynomial inequalities, which unfortunately is currently too computationally demanding in practice. 

Not all convex bodies have the Rupert property---for example, balls are not Rupert. However, since no convex polyhedron has been shown not to have the Rupert property, it has been conjectured that all convex polyhedra do:
\begin{conjecture}[Jerrard et al.~\cite{jerrard} (hesitantly), Chai et al.~\cite{chai}]
Every convex polyhedron is Rupert.
\end{conjecture}
Nevertheless, despite much computational effort, most polyhedra lack positive proof of Rupertness, and hence a contradicting conjecture has also been formulated:
\begin{conjecture}[Steininger and Yurkevich~\cite{steininger}]
The rhombicosidodecahedron is not Rupert.
\end{conjecture}

In this note, we build on the work of Steininger and Yurkevich~\cite{steininger} and utilize standard nonlinear optimization methods for finding orientations of polyhedra permitting Rupert passages. The optimization methods used here provide independent verification of the previously known results for Platonic, Archimedean, Catalan, and Johnson solids. It also establishes new results: 2 additional Catalan solids---the triakis tetrahedron and the pentagonal icositetrahedron---and 5 additional Johnson solids are also found to be Rupert.\footnote{The coordinates of the Catalan solids are taken from \url{http://www.dmccooey.com/polyhedra/Catalan.html}, and those of the Johnson solids from \url{http://www.dmccooey.com/polyhedra/Johnson.html}.} Although such optimization can show Rupertness for many polyhedra, it provides no answer as to whether all convex polyhedra are Rupert.

\section*{Preliminaries}
Although the objects of interest are polyhedra in $\mathbb{R}^3$, it is sufficient to consider their projections onto planes in $\mathbb{R}^2$ to determine whether they have the Rupert property, as described by the following theorem:
\begin{theorem}[Jerrard et al.~\cite{jerrard}]
Let $P$ be a convex body in $\mathbb{R}^3$. If there are planes $\pi_p$ and $\pi_q$ such that the projection of $P$ onto $\pi_p$ fits in the interior of the projection of $P$ onto $\pi_q$, then $P$ is Rupert.
\end{theorem}

To map polyhedra onto planes, we will use the same projection maps from $\mathbb{R}^3$ to $\mathbb{R}^2$ as Steininger and Yurkevich~\cite{steininger}. Points on the 3-dimensional sphere are parameterized by the mapping $X : [0, 2 \pi) \times [0, \pi] \to \{ x \in \mathbb{R}^3 : ||x|| = 1\}$ according to  \[
X(\theta, \phi) = (\cos \theta \sin \phi, \sin \theta \sin \phi, \cos \phi).
\]
To find a plane orthogonal to $a = X(\theta, \phi)$, we start by observing that the vector $b = (-\sin\theta, \cos\theta, 0)$ is orthogonal to $a$, and then take the cross-product of $a$ and $b$ to find the vector $c = (-\cos\theta \cos\phi, -\sin\theta \cos\phi, \sin\phi)$ that is orthogonal to both $a$ and $b$. The vectors $a, b, c$ form an orthonormal basis of $\mathbb{R}^3$, and to change basis of a vector $p$ from $x, y, z$ to $a, b, c$, we multiply $p$ by the change-of-basis matrix
\[
\left(
\begin{array}{c}
a\\
b\\
c
\end{array}
\right)
=
\left(
\begin{array}{ccc}
\cos \theta \sin \phi & \sin \theta \sin \phi & \cos \phi\\
-\sin \theta & \cos \theta & 0\\
-\cos\theta \cos\phi & -\sin\theta \cos\phi & \sin\phi
\end{array}
\right).
\]
Projections of points in 3-dimensional space onto a plane orthogonal to $X(\theta, \phi)$ can then be given by $M_{\theta, \phi} : \mathbb{R}^3 \to \mathbb{R}^2$, where
\begin{equation}\label{m-eq}
M_{\theta, \phi} = \left(
\begin{array}{ccc}
-\sin \theta & \cos \theta & 0\\
-\cos \theta \cos \phi & -\sin \theta \cos \phi & \sin \phi
\end{array}
\right),
\end{equation}
i.e., the last two rows of the change-of-basis matrix. A rotation in the plane by the angle $\alpha$ is performed by the mapping $R_\alpha : \mathbb{R}^2 \to \mathbb{R}^2$, where
\[
R_{\alpha} = \left(
\begin{array}{cc}
\cos \alpha & -\sin \alpha\\
\sin \alpha & \cos \alpha\\
\end{array}
\right),
\]
and a translation in the plane by $u$ and $v$ (in the $x$ and $y$ coordinates of the plane, respectively) is given by $T_{u,v}: \mathbb{R}^2 \to \mathbb{R}^2$, where
\[
T_{u,v}(x,y) = (x+u, y+v).
\]
The mappings $M_{\theta,\phi}$, $R_\alpha$, and $T_{x,y}$ all act on sets of points elementwise.

Denoting the interior of a set $X$ by $\mathrm{int}(X)$, an alternative to Theorem 1 can now be formulated as follows:
\begin{theorem}[Steininger and Yurkevich~\cite{steininger}]
A convex polyhedron $P$ is Rupert if there are angles $\theta_p, \phi_p$ and $\theta_q, \phi_q$, a rotation $\alpha$, and a translation $(u,v)$ such that
\begin{equation*}
\left(T_{u,v} \circ R_\alpha \circ M_{\theta_p, \phi_p}\right) ( P ) \subset \mathrm{int}\left( M_{\theta_q, \phi_q} (P) \right).
\end{equation*}
\end{theorem}
Both sides of the subset inclusion in Theorem 2 can be represented by convex polygons. This means that, given angles $\theta_p, \phi_p, \theta_q, \phi_q$, if we can find a placement (allowing for rotations and translations in the plane) of the polygon $M_{\theta_p, \phi_p}(P)$ inside the polygon $M_{\theta_q, \phi_q}(P)$, then the polyhedron $P$ is Rupert. We thus turn to finding placements of polygons.

\section*{Placement of a convex polygon inside another}\label{sec-similar}
Chazelle~\cite{chazelle} devised an algorithm for determining whether a polygon fits inside another. Since we are interested in finding good lower bounds of the Nieuwland numbers of polyhedra, we want to find as large scale factors as possible such that the polygon still fits. Steininger and Yurkevich~\cite{steininger} used Chazelle's method combined with a binary search to improve on the bounds found by their algorithm.

An alternative method that directly provides the largest scale factor for a given polygon (i.e., for a given orientation of the polyhedron) is described by Agarwal et al.~\cite{agarwal}, who utilized the nice property that in $\mathbb{R}^2$, a rotation by an angle $\alpha$ about the origin and a scaling by a factor $\rho$ can be parameterized by the Euclidean coordinates $s$ and $t$ using the polar substitution $s = \rho \cos \alpha$ and $t = \rho \sin \alpha$. They proceeded to show how the largest similar placement of a polygon $P$ inside a convex polygon $Q$ can be easily computed, given a procedure for computing the convex hull, as follows.

Let $P$ be represented by its vertices $(x_i, y_i)$, $i = 1,\ldots,m$, and $Q$ be represented by its halfplane equations $a_j x + b_j y \leq 1$, $j = 1, \ldots,n$.  A vertex $(x_i,y_i)$ of $P$ rotated by the angle $\alpha$ about the origin, scaled by the factor $\rho$, and translated by $(u,v)$ is given by
\[
(T_{u,v} \circ \rho R_\alpha)(x_i,y_i) = (x_i \rho \cos \alpha - y_i \rho \sin \alpha + u, x_i \rho \sin \alpha + y_i \rho \cos \alpha + v).
\]
Here, $\rho \cos \alpha$ and $\rho \sin \alpha$ can be equated with the Euclidean coordinates $s$ and $t$, respectively. Then
\[
(T_{u,v} \circ \rho R_\alpha)(x_i,y_i) = (x_i s  - y_i t + u, x_i t + y_i s + v).
\]
This point lies within $Q$ if 
\[
a_j(x_i s - y_i t + u) + b_j (x_i t + y_i s + v) \leq 1,
\]
or, equivalently,
\begin{equation}\label{hyper}
(a_j x_i + b_j y_i) s + ( b_j x_i - a_j y_i )t + a_j u + b_j v \leq 1
\end{equation}
for $j = 1,\ldots,n$. If this holds for all points $(x_i, y_i)$, $i=1,\ldots,m$, of $P$, then a similar copy of $P$ rotated by $\alpha$, scaled by $\rho$, and translated by $(u,v)$ is contained in $Q$. Thus, all similar placements of $P$ in $Q$ can be represented by a 4-dimensional convex polytope $C(P,Q)$ defined by the $mn$ halfspaces given by~\eqref{hyper} for $i=1,\ldots,m$, $j=1,\ldots,n$. The largest scaling $\rho$ is obtained at the point in this polytope that maximizes $\rho^2 = s^2 + t^2$. Since the maximum of a convex function over a convex polytope is attained in a vertex, it can be found by the evaluation of the function value over all vertices of $C(P,Q)$.

To enumerate the vertices of $C(P,Q)$, we use the duality transform that maps the halfspace equation $\beta_1 s + \beta_2 t + \beta_3 u + \beta_4 v \leq 1 $ of a polytope to the vertex $(\beta_1, \beta_2, \beta_3, \beta_4)$ of its dual. We also use some procedure for computing convex hulls. The vertices of $C(P,Q)$ can then be enumerated as follows:\footnote{It should be noted that Agarwal et al.~\cite{agarwal} exploit the degeneracy of $C(P,Q)$ to come up with a more elaborate and faster algorithm.}

\begin{enumerate}
    \item Use the duality transform to map the halfspace equations of $C(P,Q)$ to vertices of its dual $D$, i.e., to the set of points \[\{ (a_j x_i + b_j y_i, b_j x_i - a_j y_i, a_j, b_j) : i=1,\ldots,m, j = 1,\ldots, n\}.\]
    
    \item Compute the convex hull of the vertices of $D$, yielding its halfspace equations.
    
    \item Use the duality transform to map the so computed halfspace equations of $D$ to the vertices of $C(P,Q)$.
\end{enumerate}
These steps yield the vertices of $C(P,Q)$, and to find the largest scaling $\rho$, one simply has to check which vertex has the largest value of $\rho^2 = s^2 + t^2$.

\section*{Optimizing for the Rupert property}
A polyhedron $P$ has the Rupert property if there are angles $\theta_p, \phi_p, \theta_q, \phi_q$ such that $\rho^2 = s^2 + t^2 > 1$ for some point in $C\left(M_{\theta_p, \phi_p}(P), M_{\theta_q, \phi_q}(P) \right)$. The scaling $\rho$ is then a lower bound of the Nieuwland number of the polyhedron. Thus, we want to find angles for which the function
\begin{equation*}
f(\theta_p, \phi_p, \theta_q, \phi_q) = \max\{ s^2 + t^2 : (s,t,u,v) \in C\left(M_{\theta_p, \phi_p}(P), M_{\theta_q, \phi_q}(P) \right) \}
\end{equation*}
evaluates to a value greater than $1$. To evaluate $f(\theta_p, \phi_p, \theta_q, \phi_q)$, we compute the polygons $M_{\theta_p, \phi_p}(P)$ and $M_{\theta_q, \phi_q}(P)$ using the mapping~\eqref{m-eq}, enumerate the vertices of $C\left(M_{\theta_p, \phi_p}(P), M_{\theta_q, \phi_q}(P) \right)$ using the procedure presented as steps 1--3 in the previous section, and check which vertex has the largest value of $\rho^2 = s^2 + t^2$.

So we can evaluate $f$, and want to solve the optimization problem
\begin{equation*}
\begin{aligned}
& \underset{\theta_p, \phi_p, \theta_q, \phi_q}{\text{maximize}}
& & f(\theta_p, \phi_p, \theta_q, \phi_q)
\end{aligned}\hspace{-0.3cm}.
\end{equation*}
For brevity, we abbreviate the angles by $x$ so that $x = (\theta_p, \phi_p, \theta_q, \phi_q)$ and let $g(x) = -f(\theta_p, \phi_p, \theta_q, \phi_q)$ to turn the problem of maximizing $f$ into the problem of minimizing $g$, since minimization is the standard in many optimization solvers.

A simple optimization method that can be used to minimize $g$ is gradient descent. This method requires the gradient $\nabla g$ of $g$, which can be approximated by finite differences:
$
\nabla g(x) \approx 
\frac{1}{2h} \left(
\begin{array}{c}
g(x + h e_i) - g(x - he_i)
\end{array}
\right)_{i=1,\ldots,4},
$
where $e_i$ is the $i$th standard basis vector and $h$ is some small positive number. Now, given a point $x_n$, we can find a better point by using the update 
$
x_{n+1} = x_n - \gamma_{n+1} \nabla g(x_n),
$
where $\gamma_{n+1} > 0$ is selected so that $g(x_{n+1}) < g(x_n)$. The updating is repeated until $\lVert \nabla g(x_{n+1}) \rVert$ falls below some tolerance.

Since the objective function $g$ is nonconvex, the optimization is dependent on the starting point. To find solutions with $\rho > 1$, it is sometimes necessary to run optimizations from various starting points. To this end, the intervals $[0,2\pi)$ and $[0,\pi]$ were discretized into $k$ points each, and the optimization was run for all $k^4$ combinations of starting points, with $\theta_p$ and $\theta_q$ initialized with points from $[0,2\pi)$, and $\phi_p$ and $\phi_q$ initialized with points from $[0,\pi]$.

While gradient descent can be used to show the Rupert property of new polyhedra (the pentagonal icositetrahedron was found to be Rupert using $k=5$), there are better optimization methods. For example, one can leverage already existing implementations of more advanced optimization methods that are available in the Python package SciPy~\cite{virtanen}. Two methods that perform well on the problems under consideration are sequential quadratic programming (SLSQP in SciPy)~\cite{kraft} and Nelder-Mead~\cite{nelder}.

Sequential quadratic programming is a continuous optimization method that solves a quadratic model of the problem at hand to find a search direction, and then computes the updated variables $x_{n+1}$ by taking a step along the search direction. For an unconstrained problem like the one considered here, it amounts to using the update 
$
x_{n+1} = x_n - \gamma_{n+1} B_n^{-1} \nabla g(x_n),
$
where $B_n$ is an approximation of the Hessian of $g$, which is updated after each iteration. This update is actually that of a quasi-Newton method, which is simpler than sequential quadratic programming, and for which SciPy also provides an implementation, but SLSQP was found to be faster in this case.

The Nelder-Mead method is a derivative-free optimization method, meaning that it does not make use of derivatives to find new iterates for the variables, but only uses function value comparisons. Applied to an $n$-dimensional problem, the method maintains a set of $n+1$ points, and repeatedly tries to switch the point with the highest value of $g$ to a new point with a lower value, utilizing what is known about the other $n$ points. If no such point is found, it moves all other points towards the one with the lowest value of $g$. When the difference between the $n+1$ points in terms of variable values or function values falls below some tolerance, the method terminates and returns the point with the lowest function value.

\begin{figure}[htb]
  \centering
    \begin{subfigure}[t]{0.32\textwidth}
     \centering
      \includegraphics[width=\linewidth]{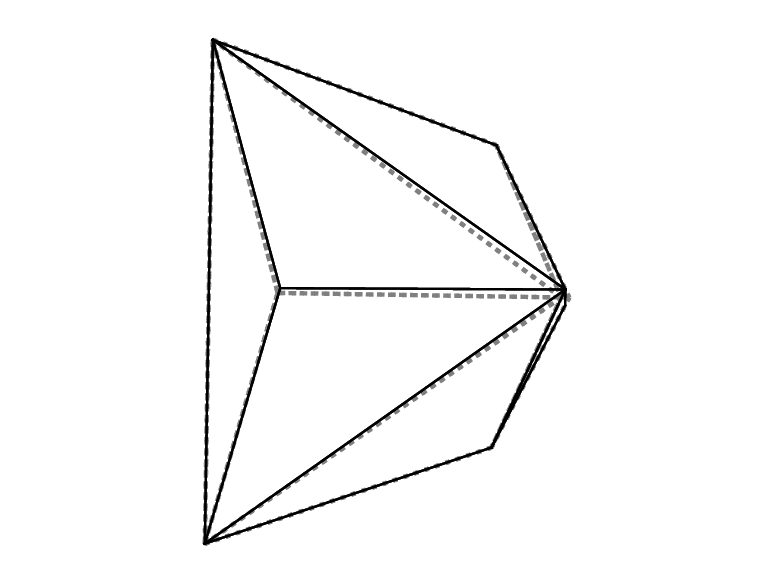}
     \caption{Triakis tetrahedron.}
    \end{subfigure}
    \begin{subfigure}[t]{0.32\textwidth}
     \centering
      \includegraphics[width=\linewidth]{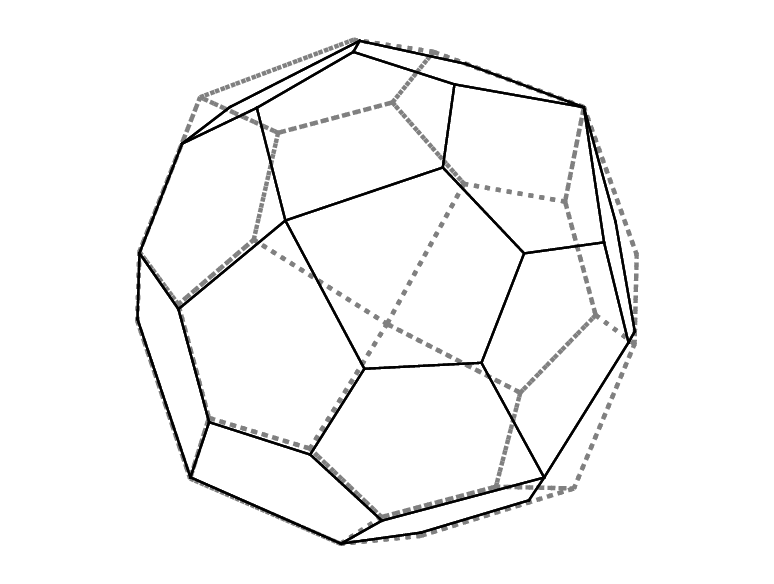}
     \caption{Pentagonal icositetrahedron.}
    \end{subfigure}
    \begin{subfigure}[t]{0.32\textwidth}
     \centering
     \includegraphics[width=\linewidth]{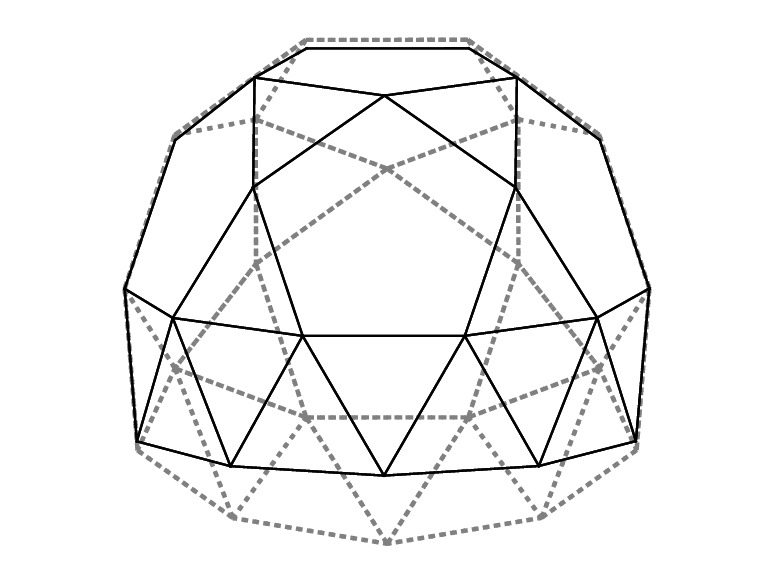} 
     \caption{J25.}
    \end{subfigure}
  \caption{New Rupert passages. Black solid lines represent polyhedra passing through those represented by gray dotted lines.}
  \label{fig-new}
\end{figure}

The optimization algorithms of SciPy were used to minimize $g$ (the source code is available at \url{https://github.com/albfre/rupert}). With $k=5$ and using Nelder-Mead with the default termination tolerance of $10^{-4}$, solutions for all Platonic, Archimedean, Catalan, and Johnson solids previously known to be Rupert were found in seconds, in most cases from one of the first few starting points. Both SLSQP and Nelder-Mead were also able to show new Rupertness results for Catalan solids (the triakis tetrahedron, the pentagonal icositetrahedron) and Johnson solids (J25, J45, J47, J71, J76). The largest values of $\rho$, i.e., the largest lower bounds for the Nieuwland numbers, were found using Nelder-Mead. Parameters for the solutions are given in Table~\ref{tab1}, and three solutions are illustrated in Figure~\ref{fig-new}.

\vspace{-0.05cm}

\begin{small}
\begin{center}
\begin{table}
\scalebox{0.86}{
\begin{tabular}{l | r r r r r r r}
 & TT & PI & J25 & J45 & J47 & J71 & J76\\
\hline
$\theta_p$ & 6.283130 & 0.466029 & 3.442081 & 3.148897 & 3.442482 &0.789632 &3.318120 \\
$\phi_p$   & 0.817234 & 1.467669 & 1.761334 & 1.568287 & 1.767543 & 1.571310 & 2.018106 \\
$\theta_q$ & 1.548107 & 2.330160 &  1.569750 & 4.719405 & 3.452852 & 3.417834 & 4.723801\\
$\phi_q$   &2.356150 & 3.026787 & 1.028319 & 2.180163 & 2.012345 & 2.726874 & 1.241903\\ 
$\alpha$   &6.2671031 &2.3256487 &0.003132 & 0.004033 & 0.001367 & 2.444476 & 5.160115 \\
$u$        &0.0001408 &  0.0006199 & 0.001326 & -0.001736 & 0.000580 & 0.004566 & 0.000365\\ 
$v$        &-0.0000022 & 0.0028453 & -0.054143 & 0.077420 & -0.033090 & -0.003943 & 0.010338\\
$\rho$     &1.000004 & 1.000436 & 1.000089 & 1.000009 & 1.000080 & 1.000598 &1.000269\\ 
\end{tabular}
}
\vspace{0.3cm}
\caption{Parameter values for new Rupert passages. TT = triakis tetrahedron, PI = pentagonal icositetrahedron. The angles $\theta_p, \phi_p, \theta_q, \phi_q$, and $\alpha$ and the translation $(u,v)$ are as in  Theorem 2, and $\rho$ is a lower bound of the Nieuwland number. (The extra digits for $\alpha, u$, and $v$ are required for TT and PI.)}
\label{tab1}
\end{table}
\end{center}
\end{small}

\vspace{-1cm}
The results were verified by explicit calculation of the two polygons $M_{\theta_p, \phi_p}(P)$ and $M_{\theta_q, \phi_q}(P)$ followed by a check that the former is contained within the latter. To determine whether a polygon with vertices $p_i =  (x_i, y_i)$, $i = 1,\ldots, m$, is contained in a convex polygon with vertices $q_j = (a_j, b_j)$, $j=1,\ldots,n$, ordered counter-clockwise, one can verify that 
\[
\mathrm{det}(q_j - p_i, q_{j+1} - p_i) =  
(a_j - x_i) (b_{j+1} - y_i) - (b_j - y_i) (a_{j+1} - x_i) > 0
\]
for all $i=1,\ldots,m$ and $j=1,\ldots,n$, where $q_{n+1}$ is defined to equal $q_1$.

Even for $k=21$ with tolerance $10^{-4}$, which took a couple of days to run for each polyhedron, and for $k=5$ with tolerance $10^{-10}$, no solution was found for any of the remaining Archimedean solids (the snub cube, the rhombicosidodecahedron, the snub dodecahedron), Catalan solids (the deltoidal hexecontahedron, the pentagonal hexecontahedron), nor Johnson solids (J72, J73, J74, J75, J77).

The difficulty of finding solutions for these polyhedra, compared to the relative easiness of finding solutions for the other ones, can make one inclined towards a belief in Conjecture 2---that there are polyhedra that are not Rupert. On the other hand, the lower bound for the Nieuwland number of the triakis tetrahedron is quite small ($1 + 4\times 10^{-6}$), and nothing says that the numbers for other polyhedra can't be even smaller. Current algorithms might be hindered by their limited numerical precision.


\end{document}